\documentclass[12pt]{article}
\usepackage{multicol}
\usepackage{amsmath, amssymb, amscd, enumerate, amsfonts, latexsym, amsxtra, 
amsthm, bm}
\usepackage{indentfirst}
\newtheorem{The}{Theorem}[section]
\newtheorem{Prop}[The]{Proposition}
\newtheorem{Lem}[The]{Lemma}

\begin{document}
\centerline{\Large Projective modules for the subalgebra of degree 0 } \vspace{3mm}
\centerline{\Large in a  finite-dimensional hyperalgebra  of type $A_1$} \vspace{7mm}
\centerline{Yutaka Yoshii 
\footnote{ E-mail address: yutaka.yoshii.6174@vc.ibaraki.ac.jp}}  \vspace{5mm}
\centerline{College of Education,   
Ibaraki University,}
\centerline{2-1-1 Bunkyo, Mito, Ibaraki, 310-8512, Japan}
\begin{abstract}
We describe the structure of projective indecomposable modules for 
the subalgebra consisting of the elements of degree 0 in the hyperalgebra of the $r$-th Frobenius kernel for the 
algebraic group ${\rm SL}_2(k)$, using the primitive idempotents which were constructed before by the author.   
\end{abstract}
\renewcommand{\thefootnote}{\fnsymbol{footnote}}

\section{Introduction}
Let $k$ be an algebraically closed field of characteristic $p>0$. Let $G$ be a connected, simply connected and semisimple algebraic group over $k$ which is split over the finite field $\mathbb{F}_p$ of order $p$. 

The representation theory of $G$ is closely related to that of   
the $r$-th Frobenius kernel $G_r$. 
Since the representation theory of $G$ can be identified with the locally finite 
representation theory of the corresponding 
(infinite-dimensional) $k$-algebra $\mathcal{U}$ which is called the hyperalgebra of $G$, and since the representation theory of 
$G_r$ can be identified with that of the corresponding finite-dimensional 
hyperalgebra $\mathcal{U}_r$, it is important to study the structure of 
projective indecomposable modules (PIMs) for  $\mathcal{U}_r$. Thus it is worthwhile constructing 
primitive idempotents in $\mathcal{U}_r$.    
Unfortunately, it seems that explicit description of primitive idempotents in $\mathcal{U}_r$ has not  been known for general $G$. If $G$ is of type $A_1$ (i.e. $G={\rm SL}(2,k)$), 
the explicit description is given in Seligman's paper \cite{seligman03} for $r=1$, and  
 in  author's paper  \cite{yoshii17}  for general $r$.    

In this paper, using the primitive idempotents in $\mathcal{U}_r$ given in the author's paper,  we shall study the projective indecomposable modules for the subalgebra  $\mathcal{A}_r$ consisting of the elements of
degree 0  in $\mathcal{U}_r$. 
More concretely, since any idempotent in $\mathcal{U}_r$ lies in $\mathcal{A}_r$, 
we can describe the structure of projective indecomposable $\mathcal{A}_r$-modules 
by giving that of the $\mathcal{A}_r$-modules generated by the primitive 
idempotents in $\mathcal{U}_r$. 
Although the argument is not so difficult, 
the structure of these modules can be completely determined. This result  enables 
us to see the primitivity of the idempotents without knowing  dimensions of the simple 
$\mathcal{U}_r$-modules. It is also expected that the $\mathcal{A}_r$-modules 
will be useful to study the structure of projective indecomposable 
$\mathcal{U}_r$-modules. 

The main results will be given in Section 3. First we construct a basis of the 
$\mathcal{A}_r$-module generated by a primitive idempotent, using a method 
which generalizes the one to construct the idempotents in \cite{yoshii17}. Then we describe the radical 
and socle series of the $\mathcal{A}_r$-modules. It turns out that each projective 
indecomposable $\mathcal{A}_r$-module is rigid and that each block algebra of 
$\mathcal{A}_r$, which corresponds to a primitive idempotent,  is symmetric. 

\section{Preliminaries}
Assume $G={\rm SL}_2(k)$ in the rest of this paper. Let 
$$X= \left( \begin{array}{cc}
         0 & 1 \\
         0 & 0 \end{array} \right),\ \ \ 
Y= \left( \begin{array}{cc}
         0 & 0 \\
         1 & 0 \end{array} \right),\ \ \ 
H= \left( \begin{array}{cc}
         1 & 0 \\
         0 & {-1} \end{array} \right)$$
be the standard basis in the simple complex Lie algebra $\mathfrak{g}_{\mathbb{C}}=\mathfrak{sl}_2(\mathbb{C})$. We define a subring 
$\mathcal{U}_{\mathbb{Z}}$ of the universal enveloping algebra $\mathcal{U}_{\mathbb{C}}$ of  $\mathfrak{g}_{\mathbb{C}}$  
generated by $X^{(m)}=X^m/m!$ and  $Y^{(m)}=Y^m/m!$ with $m \in \mathbb{Z}_{\geq 0}$. 
Set 
$${H+c \choose m }= \dfrac{(H+c)(H+c-1) \cdots (H+c-m+1)}{m!}$$
for $c \in \mathbb{Z}$ and $m \in \mathbb{Z}_{\geq 0}$. The  elements 
$$Y^{(m)} {H \choose n} X^{(m')}$$
with $m, m', n \in \mathbb{Z}_{\geq 0}$ form a $\mathbb{Z}$-basis of $\mathcal{U}_{\mathbb{Z}}$. The $k$-algebra 
$\mathcal{U}_{\mathbb{Z}} \otimes_{\mathbb{Z}} k$ is denoted by $\mathcal{U}$ or 
${\rm Dist}(G)$ which is  called 
the hyperalgebra of $G$. We  use the same notation for the images in $\mathcal{U}$ 
of the elements 
in $\mathcal{U}_{\mathbb{Z}}$. 

Let ${\rm Fr}: \mathcal{U} \rightarrow \mathcal{U}$ be 
 the $k$-algebra endomorphism 
 which is defined by 
$${\rm Fr}(X^{(m)}) = 
\left\{ \begin{array}{cl}
         {X^{(m/p)}} & {\mbox{if $p \mid m$,}} \\
         0 & {\mbox{if $p \nmid m$}} \end{array} \right. 
\ \ \ \ \ \ \mbox{and  \ }
{\rm Fr}(Y^{(m)}) = 
\left\{ \begin{array}{cl}
         {Y^{(m/p)}} & {\mbox{if $p \mid m$,}} \\
         0 & {\mbox{if $p \nmid m$}} \end{array} \right. .$$
Then we also have 
$${\rm Fr} \bigg( {H \choose m} \bigg) = 
\left\{ \begin{array}{cl}
         {{H \choose m/p}} & {\mbox{if $p \mid m$,}} \\
         0 & {\mbox{if $p \nmid m$}} \end{array} \right. .$$

Let $\mathcal{U}^0$ be 
the subalgebra of $\mathcal{U}$ generated by ${H \choose p^i}$ with 
$i \in \mathbb{Z}_{\geq 0}$. The elements 
$Y^{(m)} {H \choose n} X^{(m')}$ 
with $m, m', n \in \mathbb{Z}_{\geq 0}$ form a $k$-basis of $\mathcal{U}$. We say that an  element 
$z \in \mathcal{U}$ has degree $d$ if it is a $k$-linear combination of  the elements 
$Y^{(m)} {H \choose n} X^{(m')}$ 
with $m, m', n \in \mathbb{Z}_{\geq 0}$ and $m'-m =d$.
For a positive integer 
$r \in \mathbb{Z}_{> 0}$, let $\mathcal{U}_r$ be the subalgebra of  $\mathcal{U}$ generated by 
$X^{(p^i)}$ and $Y^{(p^i)}$ with $0 \leq i  \leq r-1$. 
This is a finite-dimensional algebra of dimension $p^{3r}$ which has 
$Y^{(m)} {H \choose n} X^{(m')}$ 
with $0 \leq m, m', n \leq p^r-1$ as a $k$-basis, 
and it can be identified with the hyperalgebra of the $r$-th Frobenius kernel 
$G_r$ of $G$. 
Let $\mathcal{U}_r^0$ be 
the subalgebra of $\mathcal{U}$ generated by ${H \choose p^i}$ with 
$0 \leq i  \leq r-1$. 

Let   ${\rm Fr}': \mathcal{U} \rightarrow \mathcal{U}$ be the $k$-linear map defined 
by 
$$Y^{(m)} {H \choose n} X^{(m')} \mapsto 
Y^{(mp)} {H \choose np} X^{(m'p)}.$$
This map is not a homomorphism of $k$-algebras, whereas its restriction to 
$\mathcal{U}_r^0$ is (for details, see \cite[\S3]{gros12} and 
\cite[\S1]{gros-kaneda11}). 
Clearly we have ${\rm Fr} \circ {\rm Fr}' = id_{\mathcal{U}}$.

Let $\mathcal{A}$ be the subalgebra of $\mathcal{U}$ 
which is generated by $\mathcal{U}^0$ and $Y^{(p^i)}X^{(p^i)}$ with 
 $i \geq 0$. This subalgebra is commutative and has the 
elements 
$Y^{(m)} {H \choose n} X^{(m)}$
with $m,  n \in \mathbb{Z}_{\geq 0}$ as a $k$-basis. For an integer 
$r >0$, set $\mathcal{A}_r = \mathcal{A} \cap \mathcal{U}_r$. This subalgebra 
is generated by $\mathcal{U}_r^0$ and $Y^{(p^i)}X^{(p^i)}$ with 
 $0 \leq i \leq r-1$, and  has the 
elements 
$Y^{(m)} {H \choose n} X^{(m)}$
with $m, n \in \{ 0,1, \dots, p^r-1\}$ as a $k$-basis. \\

For a finite-dimensional (associative) $k$-algebra $R$, let ${\rm rad}R$ be the 
largest nilpotent ideal of $R$ which is called the Jacobson radical of $R$. For a 
finite-dimensional (left) $R$-module $M$ and a positive integer $n$, 
the $R$-submodule $({\rm rad} R)^n M$ is denoted by ${\rm rad}_R^n M$ and called 
the $n$-th radical of $M$. For convenience, set ${\rm rad}_R^0 M=M$. 
If $n=1$,  the submodule ${\rm rad}_R^1 M$ is 
denoted by ${\rm rad}_R M$ and called the radical of $M$. 

In turn, for a finite-dimensional (left) $R$-module $M$ and a positive integer $n$, 
the $R$-submodule of $M$ consisting of the elements annihilated by $({\rm rad} R)^n$ is denoted 
by ${\rm soc}_R^n M$, which is called the $n$-th socle of $M$. 
For convenience, set ${\rm soc}_R^0 M=0$. 
If $n=1$,  the submodule ${\rm soc}_R^1 M$ is 
denoted by ${\rm soc}_R M$ and called the socle of $M$, which is also the 
largest semisimple $R$-submodule of $M$.

\section{PIMs for $\mathcal{A}_r$}
 To study the projective indecomposable $\mathcal{A}_r$-modules, we use  
the idempotents constructed in  author's paper \cite{yoshii17}. 

For $a \in \mathbb{Z}$, set 
$$\mu_a = {H-a-1 \choose p-1} = \sum_{i=0}^{p-1}
{-a-1 \choose p-1-i} {H \choose i} \in \mathcal{U}_1^0.$$
This is a   $\mathcal{U}_1^0$-weight vector of weight $a$ in the $\mathcal{U}_1^0$-module  
$\mathcal{U}_1^0$: $H \mu_a = a \mu_a$.
Moreover,  we have  $\mu_a = \mu_b$ if and only if $a \equiv b\ ({\rm mod}\ p)$, and 
all $\mu_a$ with $a \in \{ 0,1, \dots, p-1\}$ are pairwise orthogonal 
primitive idempotents in $\mathcal{U}_1^0$ whose sum is $1 \in \mathcal{U}_1^0$.

Suppose for a moment  that $p$ is odd. 
Set $\mathcal{S} = \{ 0, 1, \dots , (p-1)/2 \} \subset \mathbb{Z}$. 
We denote by $\mathcal{S}$ again the image of the subset 
$\mathcal{S} \subset \mathbb{Z}$ under the natural map 
$\mathbb{Z} \rightarrow \mathbb{F}_p$. 
For $\varepsilon \in \{ 0,1 \}$, $a \in \mathbb{Z}$, $j \in \mathcal{S}$ 
and $m \in \mathbb{Z}_{\geq 0}$ we define  polynomials $\varphi_{a,m}(x), \psi(x), \psi_j^{(\varepsilon)}(x) \in \mathbb{F}_p [x]$  as  
$$\varphi_{a,0}(x)= 1,$$
$$\varphi_{a,m}(x)= \prod_{i=0}^{m-1} \big(x-i(i+a+1)\big)$$
if $m \neq 0$, 
$$\psi(x)=\prod_{i \in \mathbb{F}_p} (x- i^2) = x \prod_{i \in \mathcal{S}- \{ 0\}} 
(x- i^2)^2, $$
$$\psi_j^{(1)}(x) = \psi(x) / (x- j^2),$$
$$\psi_0^{(0)}(x)= \prod_{i \in \mathbb{F}_p^{\times}} (x-i^2) = 
\prod_{i \in \mathcal{S}-\{ 0\}} (x-i^2)^2$$
and 
$$\psi_{s}^{(0)}(x) = 2x(x+s^2) \prod_{i \in \mathbb{F}_p^{\times} - \{ s, p-s \} } (x-i^2)
 = 2x(x+s^2) \prod_{i \in \mathcal{S} - \{0, s \} } (x-i^2)^2 
$$
if $s \neq 0$. Clearly $\psi_0^{(0)}(x)=\psi_0^{(1)}(x)$, and we have  
$$\psi \Big(x+\big((a+1)/2\big)^2\Big) = \varphi_{a,p}(x)$$ and  
$$\varphi_{a,p} (\mu_a YX) = \varphi_{-a, p}(\mu_a X Y) =0$$
(see \cite[Lemma 4.3]{yoshii17}).  Set $\mathcal{P} = \{0,1, \dots, p-1\} \times \mathcal{S}$
and 
$$B^{(\varepsilon)}(a, j) = \psi_j^{(\varepsilon)} \Big(\mu_a YX +\big((a+1)/2\big)^2\Big) \cdot \mu_a$$
for $\varepsilon \in \{ 0,1 \}$ and $(a,j) \in \mathcal{P}$.  
This element also can be written as 
$$B^{(\varepsilon)}(a, j) = \psi_j^{(\varepsilon)} \Big(\mu_a XY +\big((a-1)/2\big)^2\Big) \cdot \mu_a.$$ Note also that $B^{(0)}(a,0)=B^{(1)}(a,0)$ for any 
$a \in \{0,1, \dots, p-1 \}$. 

In turn, suppose that  $p=2$. Then we consider the set 
$$\mathcal{P}=\{ (0,1/2), (1,0), (1,1)\} \subset 
\{0,1 \} \times (1/2)\mathbb{Z}$$ instead of $\mathcal{P}=\{ 0,1, \dots, p-1\} \times \mathcal{S}$ when $p$ is odd, and define 
$$B^{(0)}(0,1/2)= \mu_0,\ \ \ B^{(1)}(0,1/2)= \mu_0YX= \mu_0 XY, $$
$$B^{(0)}(1,0)=B^{(1)}(1,0)= \mu_1YX = \mu_1 XY + \mu_1,$$
$$B^{(0)}(1,1)=  B^{(1)}(1,1)= \mu_1YX + \mu_1=\mu_1XY.$$

For any prime number $p$ and a pair $(a,j) \in \mathcal{P}$, set $E(a,j) = B^{(0)}(a,j)$. The elements $E(a,j)$ with $(a,j) \in \mathcal{P}$ are 
 pairwise orthogonal idempotents in $\mathcal{U}_1$ whose sum is   the unity $1 \in \mathcal{U}_1$ 
(see \cite[Proposition 4.5]{yoshii17}).    \\

To construct idempotents in $\mathcal{A}_r$, we make 
some preparation.

First we shall define $n^{(\varepsilon)}(a,j)$ for each 
$\varepsilon \in \{0,1\}$ and a pair $(a,j)$ in $\mathbb{Z} \times \mathcal{S}$ 
(when $p$ is odd) or  $\mathcal{P}$ (when $p=2$) as follows: 
if $p$ is odd, $n^{(\varepsilon)}(a,j)$ is the largest 
non-negative integer $n$ satisfying 
$$\varphi_{a,n} (x) \mid \psi^{(\varepsilon)}_j \Big(x + \big((a+1)/2\big)^2\Big)$$ 
for $(a,j) \in \mathbb{Z} \times \mathcal{S}$, and if $p=2$, we set 
$$n^{(0)}(0, 1/2) = 0,\ \ n^{(0)}(1,0)=1,\ \ n^{(0)}(1,1)=0,$$
$$n^{(1)}(0, 1/2) = 1,\ \ n^{(1)}(1,0)=1,\ \ n^{(1)}(1,1)=0.$$

We consider the following four conditions for each pair $(a,j) \in \mathcal{P}$: \\ \\
(A) $a$ is even and $(p-a+1)/2 \leq j \leq (p-1)/2$,\\
(B) $a$ is even and $0 \leq j \leq (p-a-1)/2$,\\
(C) $a$ is odd and $0 \leq j \leq (a-1)/2$,\\
(D) $a$ is odd and $(a+1)/2 \leq j \leq (p-1)/2$.\\ \\
Note that if $p=2$, the pairs $(0,1/2)$, $(1,0)$ and $(1,1)$ in $\mathcal{P}$ satisfy (B), (C) and (D) 
respectively.  

\begin{Lem}
Let $(a,j) \in \mathcal{P}$. 
Then the following holds. \\ \\
{\rm (i)} $n^{(0)}(a,j)= (p-a-1)/2 +j\ \mbox{ and }
\ n^{(1)}(a,j)= (3p-a-1)/2 -j$ under {\rm (A)}, \\
{\rm (ii)} $n^{(0)}(a,j)= (p-a-1)/2 -j\ \mbox{ and }
\ n^{(1)}(a,j)= (p-a-1)/2 +j$ under   {\rm (B)}, \\
{\rm (iii)} $n^{(0)}(a,j)= (2p-a-1)/2 -j \ 
\mbox{ and } \ 
n^{(1)}(a,j)= (2p-a-1)/2 +j$ under {\rm (C)}, \\
{\rm (iv)} $n^{(0)}(a,j)= j- (a+1)/2\ \mbox{ and } 
\ n^{(1)}(a,j)= (2p-a-1)/2 -j$ under  {\rm (D)}.
\end{Lem}
\noindent {\itshape Proof.} It is clear when $p=2$. Suppose that $p$ is odd. 
Then $n^{(0)}(a,j)$ is determined in \cite[Lemma 4.6]{yoshii17}, and so 
we only have to determine 
$n^{(1)}(a,j)$. We have $n^{(0)}(a,0)=n^{(1)}(a,0)$ since 
$\psi_0^{(0)}(x)=\psi_0^{(1)}(x)$, and so the lemma holds for $j =0$. Suppose that 
$j \neq 0$. Then the definition of 
$n^{(1)}(a,j)$ implies that it is the second smallest non-negative integer $n$ 
satisfying  $x-n(n+a+1) = x+\big( (a+1)/2\big)^2 -j^2$ in $\mathbb{F}_p[x]$, and hence 
$-n(n+a+1) = \big( (a+1)/2\big)^2 -j^2$ in $\mathbb{F}_p$. Thus we obtain the result for 
$n^{(1)}(a,j)$ (see the second to fifth paragraphs in the proof of  \cite[Lemma 4.6]{yoshii17}). 
$\square$ \\ 

If $p$ is odd, we define an integer $\tilde{n}^{(\varepsilon)}(a,j)$ for 
$\varepsilon \in \{0,1 \}$, $a \in \mathbb{Z}$ and $j \in \mathcal{S}$  
as $\tilde{n}^{(\varepsilon)}(a,j) = n^{(\varepsilon)}(-a,j)$. If $p=2$, set  
$$\tilde{n}^{(0)}(0,1/2)=0,\ \ \tilde{n}^{(0)}(1,0) = 0,\ \ \tilde{n}^{(0)}(1,1)=1,$$ 
$$\tilde{n}^{(1)}(0,1/2)=1,\ \ \tilde{n}^{(1)}(1,0) = 0 ,\ \ \tilde{n}^{(1)}(1,1)=1.$$ 
It is easy to see that  Lemma 3.1 implies the following (see \cite[Lemma 4.7]{yoshii17}). 
\\

\begin{Lem}
Let $(a,j) \in \mathcal{P}$. Then the following holds. \\ \\
{\rm (i)} 
$\tilde{n}^{(0)}(a,j)= (-p+a-1)/2 +j\ \mbox{ and }
\ \tilde{n}^{(1)}(a,j)= (p+a-1)/2 -j.$  under {\rm (A)}, \\
{\rm (ii)} 
$\tilde{n}^{(0)}(a,j)= (p+a-1)/2 -j\ \mbox{ and }
\ \tilde{n}^{(1)}(a,j)= (p+a-1)/2 +j$ under  {\rm (B)}, \\
{\rm (iii)} 
$\tilde{n}^{(0)}(a,j)= (a-1)/2 -j\ \mbox{ and }
\ \tilde{n}^{(1)}(a,j)= (a-1)/2 +j$ under {\rm (C)},\\ 
{\rm (iv)} 
$\tilde{n}^{(0)}(a,j)= (a-1)/2+j \ \mbox{ and }
\ \tilde{n}^{(1)}(a,j)= (2p+a-1)/2 -j$ under   {\rm (D)}.
\end{Lem}

The following lemma is a generalization of \cite[Lemma 4.8]{yoshii17}.

\begin{Lem}
Let $(a,j) \in \mathcal{P}$. Then the element $B^{(\varepsilon)}(a,j)$ with 
$\varepsilon \in \{0,1 \}$ can be written as
$$B^{(\varepsilon)}(a,j)= \mu_a \sum_{m=n^{(\varepsilon)}(a,j)}^{p-1} 
c^{(\varepsilon)}_m Y^m X^m = 
\mu_a \sum_{m=\tilde{n}^{(\varepsilon)}(a,j)}^{p-1} \tilde{c}^{(\varepsilon)}_m X^m Y^m$$
for some $c^{(\varepsilon)}_m, \tilde{c}^{(\varepsilon)}_m \in \mathbb{F}_p$ with 
$c^{(\varepsilon)}_{n^{(\varepsilon)}(a,j)} \neq 0$ and 
$\tilde{c}^{(\varepsilon)}_{\tilde{n}^{(\varepsilon)}(a,j)} \neq 0$. 
\end{Lem}
\noindent {\itshape Proof.} The equalities for $p=2$ are clear by the definition of 
$B^{(\varepsilon)}(a,j)$. Suppose that  $p$ is odd.  Then the equalities for 
$E(a,j)=B^{(0)}(a,j)$ are  proved in \cite[Lemma 4.8]{yoshii17}, 
using the integers $n^{(0)}(a,j)$, 
$\tilde{n}^{(0)}(a,j)$ and the polynomial 
$\psi_j^{(0)}(x)$. The proof for $B^{(1)}(a,j)$ is similar,  using $n^{(1)}(a,j)$, 
$\tilde{n}^{(1)}(a,j)$ and $\psi_j^{(1)}(x)$.  $\square$  \\

Now we shall construct idempotents in $\mathcal{A}_r$ for $r \in \mathbb{Z}_{>0}$. 
First of all, we give the primitive idempotents in $\mathcal{U}^0_r$. For an integer 
$a \in \mathbb{Z}$, set 
$$\mu_a^{(r)} = {H-a-1 \choose p^r-1} \in \mathcal{U}^0_{r}.$$
This is a   $\mathcal{U}_r^0$-weight vector of weight $a$ in the $\mathcal{U}_r^0$-module  
$\mathcal{U}_r^0$:  
${H \choose n} \mu_a = {a \choose n} \mu_a$ for any $n \in \{0,1,\dots, p^r-1 \}$.
Moreover,  we have  $\mu^{(r)}_a = \mu^{(r)}_b$ if and only if $a \equiv b\ ({\rm mod}\ p^r)$, and 
all $\mu^{(r)}_a$ with $a \in \{ 0,1, \dots, p^r-1\}$ are pairwise orthogonal 
primitive idempotents in $\mathcal{U}_r^0$ whose sum is $1 \in \mathcal{U}_r^0$. 
For details, see \cite[\S4.7]{gros-kaneda15}.

If $(a,j) \in \mathcal{P}$ satisfies {\rm (A)} or 
{\rm (C)}, then set 
$s(a,j)=(p-a+1)/2$ if $p$ is odd and $a$ is even,  $s(a,j)=(p-a)/2$ if both $p$ and $a$ are  odd, 
and  $s(a,j)=1$ if $p=2$. 

For $\varepsilon \in \{0,1 \}$ and $(a,j) \in \mathcal{P}$,  we write 
$$B^{(\varepsilon)}(a,j)= \mu_a \sum_{m=n^{(\varepsilon)}(a,j)}^{p-1} 
c^{(\varepsilon)}_m Y^m X^m = 
\mu_a \sum_{m=\tilde{n}^{(\varepsilon)}(a,j)}^{p-1} \tilde{c}^{(\varepsilon)}_m X^m Y^m$$
following Lemma 3.3. 
Using this notation 
we define $Z^{(\varepsilon)} \big(z; (a,j) \big)$ for $z \in \mathcal{U}$ as 
$$Z^{(\varepsilon)} \big(z; (a,j) \big)
= \mu_a \sum_{m=n^{(\varepsilon)}(a,j)}^{p-1} 
c^{(\varepsilon)}_m Y^m X^{m-s(a,j)} {\rm Fr'}(z)X^{s(a,j)}$$
if $(a,j)$ satisfies (A) or (C), and
$$Z^{(\varepsilon)} \big(z; (a,j) \big)
={\rm Fr'}(z) B^{(\varepsilon)}(a,j)$$
if $(a,j)$ satisfies (B) or (D). 
Clearly the map $Z^{(\varepsilon)} \big(-; (a,j) \big) : \mathcal{U} \rightarrow \mathcal{U}, \ 
z \mapsto Z^{(\varepsilon)} \big(z; (a,j) \big)$ is  $k$-linear. 

We introduce the following two lemmas to prove the main result. 

\begin{Lem}
For a pair $(a,j) \in \mathcal{P}$ and a nonzero element $z \in \mathcal{U}$, there is  
a nonzero element $z' \in \mathcal{U}$ which is independent of 
$\varepsilon \in \{0,1 \}$ such that $$Z^{(\varepsilon)} \big(z; (a,j) \big)
={\rm Fr'}(z') B^{(\varepsilon)}(a,j)=B^{(\varepsilon)}(a,j){\rm Fr'}(z').$$
Moreover, if $z \in \mathcal{A}$, then $z'$ and $Z^{(\varepsilon)} \big(z; (a,j) \big)$ also lie in $\mathcal{A}$.  
\end{Lem} 
\noindent {\bf Remark.} This lemma implies the following  facts: \\ \\
(a) If $p$ is odd and $j=0$, or if $p=2$ and $a=1$, then we have 
$$Z^{(0)} \big(z; (a,j) \big) = Z^{(1)} \big(z; (a,j) \big)$$ since 
$B^{(0)}(a,j) = B^{(1)}(a,j)$. Otherwise $Z^{(0)} \big(z_0; (a,j) \big)$ and 
$Z^{(1)} \big(z_1; (a,j) \big)$ are linearly independent over $k$  for 
nonzero elements $z_0,z_1 \in \mathcal{U}$. \\ \\
(b) The $k$-linear map  $Z^{(\varepsilon)} \big(-; (a,j) \big)$ is injective.  

\begin{Lem}
Let $u$ be an element of the $k$-subalgebra of $\mathcal{U}$ generated by all $X^{(p^i)}$ and $Y^{(p^i)}$ 
with $i \in \mathbb{Z}_{> 0}$. For $(a,j) \in \mathcal{P}$, 
$\varepsilon \in \{0,1 \}$ and $z \in \mathcal{U}$, 
we have 
$$u Z^{(\varepsilon)} \big(z; (a,j) \big)= 
Z^{(\varepsilon)} \big( {\rm Fr}(u)z; (a,j) \big). $$
\end{Lem}
These lemmas are generalizations of Lemma 5.3 and Proposition 5.4 (iv) in 
\cite{yoshii17}. We can prove them similarly since Lemma 5.2 in \cite{yoshii17} can 
be applied even if $\varepsilon =1$ (note that $n^{(1)}(a,j) \geq n^{(0)}(a,j)$ and 
$\tilde{n}^{(1)}(a,j) \geq \tilde{n}^{(0)}(a,j)$ by Lemmas 3.1 and 3.2). 

For a positive integer $r$, consider an $r$-tuple 
$\big((a_i, j_i)\big)_{i=0}^{r-1} =\big( (a_0,j_0),\dots, (a_{r-1},j_{r-1}) \big)$ of  pairs 
$(a_i,j_i) \in \mathcal{P}$ $(0 \leq i \leq r-1)$. For convenience we shall write this as 
$$((a_0,  \dots, a_{r-1}),(j_0,  \dots, j_{r-1})),$$
or $({\bm a},{\bm j})$ with ${\bm a}=(a_0,  \dots, a_{r-1})$ and 
${\bm j}=(j_0,  \dots, j_{r-1})$.

For $r$-tuples $(\varepsilon_0, \dots, \varepsilon_{r-1}) \in \{ 0,1\}^r$ and  
$((a_0,  \dots, a_{r-1}),(j_0,  \dots, j_{r-1})) \in \mathcal{P}^r$,  
we define an element 
$B^{(\varepsilon_0, \dots, \varepsilon_{r-1})}((a_0,  \dots, a_{r-1}),(j_0,  \dots, j_{r-1})) 
\in \mathcal{U}$ 
as $B^{(\varepsilon_0)}(a_0,j_0) $ if $r=1$, and 
$$
Z^{(\varepsilon_0)} \Big( 
B^{(\varepsilon_1, \dots, \varepsilon_{r-1})}\big((a_1,  \dots, a_{r-1}),(j_1,  \dots, j_{r-1})\big); (a_0,j_0) \Big)
$$
if $r \geq 2$. We often denote this element by $B^{(\bm{\varepsilon} )}({\bm a},{\bm j})$ 
with $\bm{\varepsilon}=(\varepsilon_0, \dots, \varepsilon_{r-1})$, ${\bm a}=(a_0,  \dots, a_{r-1})$ and 
${\bm j}=(j_0,  \dots, j_{r-1})$. Clearly all $B^{(\varepsilon_0, \dots, \varepsilon_{r-1})}((a_0,  \dots, a_{r-1}),(j_0,  \dots, j_{r-1}))$ lie in $\mathcal{A}_r$. 

As in \cite[Proposition 5.5 (i)]{yoshii17}, for $\bm{\varepsilon}=(\varepsilon_0, \dots, \varepsilon_{r-1})$ and an $r$-tuple 
$({\bm a}, {\bm j})=\big( (a_i, j_i) \big)_{i=0}^{r-1} \in \mathcal{P}^r$, 
 the element 
$B^{(\bm{\varepsilon} )}({\bm a},{\bm j})$ 
is a $\mathcal{U}_r^0$-weight vector of 
 $\mathcal{U}_r^0$-weight $\sum_{i=0}^{r-1}b_i p^i$, where 
$$b_i = 
\left\{ \begin{array}{ll} 
a_i-p & {\mbox{if $(a_i,j_i)$ satisfies {\rm (A)} or {\rm (C)},}} \\
a_i   & {\mbox{if $(a_i,j_i)$ satisfies {\rm (B)} or {\rm (D)}}}
\end{array} \right. $$
since $\mu_{\sum_{i=0}^{r-1}b_i p^i}^{(r)} B^{(\bm{\varepsilon} )}({\bm a},{\bm j}) = 
B^{(\bm{\varepsilon} )}({\bm a},{\bm j})$.

The following proposition is used to remove duplicates from the elements 
$B^{(\bm{\varepsilon} )}({\bm a},{\bm j})$ with $\bm{\varepsilon} \in \{ 0,1\}^r$. 

\begin{Prop}
For $\bm{\varepsilon}=(\varepsilon_0, \dots, \varepsilon_{r-1}),  \tilde{\bm{\varepsilon}}=(\tilde{\varepsilon}_0, \dots, \tilde{\varepsilon}_{r-1}) \in 
\{0,1 \}^r$ and an $r$-tuple 
$({\bm a}, {\bm j})=\big( (a_i, j_i) \big)_{i=0}^{r-1} \in \mathcal{P}^r$, 
we have $B^{(\bm{\varepsilon} )}({\bm a},{\bm j})= 
B^{(\tilde{\bm{\varepsilon}} )}({\bm a},{\bm j})$  
if $\varepsilon_s = \tilde{\varepsilon}_s$ for any integer $s$ satisfying $j_s \neq 0$ 
when $p$ is odd or $a_s \neq 1$ when $p=2$. 
\end{Prop} 

\noindent {\itshape Proof.}  If $r=1$, the proposition holds since 
$B^{(0)} (a,0) = B^{(1)}(a,0)$  for any 
$a \in \{ 0,1, \dots , p-1\}$ when $p$ is odd, and since 
$B^{(0)} (1,j) = B^{(1)}(1,j)$ for any $j \in \{0,1 \}$ when $p=2$. 

Suppose that $r \geq 2$  and that 
$\varepsilon_s = \tilde{\varepsilon}_s$ if $j_s \neq 0$ when $p$ is odd, or if  
 $a_s \neq 1$ when $p=2$.   By induction, we have 
\begin{eqnarray*}
B^{(\bm{\varepsilon} )}({\bm a},{\bm j}) & = &
Z^{(\varepsilon_0)}\Big( 
B^{(\varepsilon_1, \dots, \varepsilon_{r-1})}\big((a_1,  \dots, a_{r-1}),(j_1,  \dots, j_{r-1})\big); (a_0,j_0) \Big) \\
& = & Z^{(\varepsilon_0)}\Big( 
B^{(\tilde{\varepsilon}_1, \dots, \tilde{\varepsilon}_{r-1})}\big((a_1,  \dots, a_{r-1}),(j_1,  \dots, j_{r-1})\big); (a_0,j_0) \Big). 
\end{eqnarray*}
Thus if $j_0 \neq 0$ when $p$ is odd, or if  $a_0 \neq 1$ when $p=2$, we have $\varepsilon_0 = \tilde{\varepsilon}_0$ and 
hence $B^{(\bm{\varepsilon} )}({\bm a},{\bm j})= 
B^{(\tilde{\bm{\varepsilon}} )}({\bm a},{\bm j})$. On the other hand, if 
$j_0 =0$ when $p$ is odd, or if  $a_0=1$ when $p=2$, 
we have $B^{(0)}(a_0,j_0)=B^{(1)}(a_0,j_0)$ and hence 
$B^{(\bm{\varepsilon} )}({\bm a},{\bm j})= 
B^{(\tilde{\bm{\varepsilon}} )}({\bm a},{\bm j})$ by (a) in the remark  of Lemma 3.4. $\square$  \\

For an $r$-tuple 
$({\bm a}, {\bm j})=\big( (a_i, j_i) \big)_{i=0}^{r-1} \in \mathcal{P}^r$, set 
$E({\bm a}, {\bm j})= 
B^{(0, \dots, 0)}({\bm a}, {\bm j}).$ 
The elements $E({\bm a}, {\bm j})$ are pairwise orthogonal 
idempotents in $\mathcal{U}_r$ whose sum is the unity $1 \in \mathcal{U}_r$. 
Actually it turns out that these idempotents are primitive, 
since we know the dimensions of all simple $\mathcal{U}_r$-modules (see 
 \cite[Proposition 5.5 (iii)]{yoshii17}). In this paper 
we will see the primitivity as the result in Theorem 3.11, without using them.

Let ${\bm e}_i$ denote an element of $\mathbb{Z}^r$ with  $1$ in the  $i$-th entry and 
$0$ elsewhere. 

\begin{Prop}
For 
${\bm \varepsilon} =(\varepsilon_0, \dots, \varepsilon_{r-1}) \in \{0,1 \}^r$, 
$({\bm a},{\bm j})=\big((a_i,j_i)\big)_{i=0}^{r-1} \in \mathcal{P}^{r}$ and an integer 
$s$ with $0 \leq s \leq r-1$, 
$Y^{(p^s)} X^{(p^s)}B^{({\bm \varepsilon})} ({\bm a},{\bm j})$ is equal to 
$$\Bigg( j_s^2 - \bigg( \dfrac{a_{s}+1}{2}\bigg)^2\Bigg) B^{({\bm \varepsilon})} ({\bm a},{\bm j}) +4j_s^2 
B^{({\bm \varepsilon }+{\bm e}_{s+1})} ({\bm a},{\bm j})$$
if  $\varepsilon_s =0$, and to 
$$\Bigg( j_s^2 - \bigg( \dfrac{a_{s}+1}{2}\bigg)^2\Bigg) B^{({\bm \varepsilon})} 
({\bm a},{\bm j})$$
if  $\varepsilon_s =1$. \\ 
\end{Prop}
\noindent {\bf Remark.} The coefficients $j_s^2 - \big( ({a_{s}+1})/{2}\big)^2$ and $4j_s^2$ make 
sense as elements in $\mathbb{F}_p$ even if $p=2$. Indeed, they are integers since 
$(a_s,j_s) \in \mathcal{P}=\{(0,1/2),(1,0),(1,1) \}$. \\

\noindent {\itshape Proof.} 
We use induction on $r$. 
If $r=1$, we easily see that  
$$Y X B^{(0)} (a_0,j_0) = 
\Bigg( j_0^2 - \bigg( \dfrac{a_{0}+1}{2}\bigg)^2\Bigg) B^{(0)} (a_0,j_0) +4j_0^2 
B^{(1)} (a_0,j_0)$$
and 
$$Y X B^{(1)} (a_0,j_0) = 
\Bigg( j_0^2 - \bigg( \dfrac{a_{0}+1}{2}\bigg)^2\Bigg) B^{(1)} (a_0,j_0)$$
by the definition of $B^{(\varepsilon_0)} (a_0,j_0)$ and the claim follows. 

Suppose that $r \geq 2$. By Lemma 3.4, there exists an element $z' \in \mathcal{A}$ 
which is independent of $\varepsilon_0$ such that 
$$B^{({\bm \varepsilon})} ({\bm a},{\bm j}) = {\rm Fr}'(z') B^{(\varepsilon_0)}(a_0,j_0)
=  B^{(\varepsilon_0)}(a_0,j_0) {\rm Fr}'(z').$$
This shows the desired equality for $s=0$ as in the last paragraph, 
so we may assume $s \geq 1$. 
Set 
${\bm \varepsilon}' = (\varepsilon_1,\dots, \varepsilon_{r-1}) \in \{0,1 \}^{r-1}$ 
and $({\bm a}', {\bm j}') = \big((a_i, j_i)\big)_{i=1}^{r-1} \in \mathcal{P}^{r-1}$. 
By Lemma 3.5 we have 
\begin{eqnarray*}
Y^{(p^s)} X^{(p^s)} B^{({\bm \varepsilon})}({\bm a}, {\bm j})  
& = &  Y^{(p^s)} X^{(p^s)} Z^{(\varepsilon_0)} \Big( 
B^{({\bm \varepsilon}')}({\bm a}', {\bm j}') 
; (a_0,j_0) \Big) \\
& = & Z^{(\varepsilon_0)} \Big( 
Y^{(p^{s-1})} X^{(p^{s-1})} B^{({\bm \varepsilon}')}({\bm a}', {\bm j}') 
; (a_0,j_0) \Big). 
\end{eqnarray*}
By induction, $Y^{(p^{s-1})} X^{(p^{s-1})} B^{({\bm \varepsilon}')}({\bm a}', {\bm j}')$ 
is equal to 
$$\Bigg( j_s^2 - \bigg( \dfrac{a_{s}+1}{2}\bigg)^2\Bigg) B^{({\bm \varepsilon}')} ({\bm a}',{\bm j}') +4j_s^2 
B^{({\bm \varepsilon }'+{\bm e}_{s}')} ({\bm a}',{\bm j}')$$
if $\varepsilon_s =0$, and to 
$$\Bigg( j_s^2 - \bigg( \dfrac{a_{s}+1}{2}\bigg)^2\Bigg) B^{({\bm \varepsilon}')} 
({\bm a}',{\bm j}')$$
if  $\varepsilon_s =1$,  where 
${\bm e}'_i$ denotes an element of $\mathbb{Z}^{r-1}$ with  $1$ in the  $i$-th entry and 
$0$ elsewhere.  Now the proposition follows from the linearity of the map 
$Z^{(\varepsilon_0)} \Big( -
; (a_0,j_0) \Big)$. 
$\square$ \\ \\ 

A partial order in $\{0,1 \}^{r}$ can be defined as 
$$(\varepsilon_0, \dots, \varepsilon_{r-1}) \leq 
(\tilde{\varepsilon}_0, \dots, \tilde{\varepsilon}_{r-1})\ 
\mbox{if $\varepsilon_i \leq \tilde{\varepsilon}_i$ for each $i$.}$$
For  ${\bm m} 
= (m_0, \dots, m_{r-1}), 
\tilde{\bm m}= (\tilde{m}_0, \dots, \tilde{m}_{r-1}) 
\in \mathbb{Z}^r$,  define the Hamming distance 
$d({\bm m},  \tilde{\bm m})$ of ${\bm m}$ and 
$\tilde{\bm m}$ as the number of the integers $i$ with 
$m_i \neq \tilde{m}_i$, and the Hamming weight 
$\mathcal{W}({\bm m})$ of ${\bm m}$ as the number of the integers 
$i$ with 
$m_i \neq 0$. 

For an $r$-tuple $({\bm a}, {\bm j}) =\big( (a_i,j_i)\big)_{i=0}^{r-1} \in \mathcal{P}^r$,  define a subset 
$\mathcal{X}_r({\bm a},{\bm j})$ of $\{0,1 \}^r$ as follows: 
$$\mathcal{X}_r({\bm a}, {\bm j}) = 
\{ (\varepsilon_0, \dots, \varepsilon_{r-1})\in \{0,1 \}^r\ |\ 
\varepsilon_i=0\ \mbox{whenever $j_i=0$}\} $$
if $p$ is odd, and 
$$\mathcal{X}_r({\bm a}, {\bm j}) = 
\{ (\varepsilon_0, \dots, \varepsilon_{r-1})\in \{0,1 \}^r\ |\ 
\varepsilon_i=0\ \mbox{whenever $a_i=1$}\} $$
if $p=2$.

From now on we shall fix 
$({\bm a}, {\bm j})= \big( (a_i,j_i)\big)_{i=0}^{r-1} \in \mathcal{P}^r$ unless 
otherwise stated in order to 
study the structure of  the $\mathcal{A}_r$-module 
$\mathcal{A}_r \cdot E({\bm a}, {\bm j})$. 

\begin{The} For ${\bm \varepsilon} \in \mathcal{X}_r({\bm a}, {\bm j})$, 
the elements 
$B^{({\bm \theta })} ({\bm a},{\bm j})$ with 
${\bm \theta} \in \mathcal{X}_r({\bm a},{\bm j})$ and 
${\bm \theta} \geq {\bm \varepsilon}$ form a $k$-basis of the 
$\mathcal{A}_r$-module  $\mathcal{A}_r \cdot 
B^{({\bm \varepsilon})}({\bm a},{\bm j})$. 
\end{The}
\noindent {\bf Remark.}  This theorem implies some facts: \\ \\
(a) For $({\bm a}, {\bm j}),(\tilde{\bm a}, \tilde{\bm j}) \in \mathcal{P}^r$,  
 ${\bm \varepsilon} \in \mathcal{X}_r({\bm a},{\bm j})$ and  $\tilde{\bm \varepsilon} \in \mathcal{X}_r(\tilde{\bm a},\tilde{\bm j}) $, we have 
$$B^{({\bm \varepsilon })} ({\bm a},{\bm j}) 
B^{(\tilde{\bm \varepsilon })} (\tilde{\bm a},\tilde{\bm j}) =0$$
if $({\bm a}, {\bm j}) \neq (\tilde{\bm a}, \tilde{\bm j})$, since 
$B^{({\bm \varepsilon })} ({\bm a},{\bm j}) \in \mathcal{A}_r \cdot E({\bm a},{\bm j})$ 
and 
$B^{(\tilde{\bm \varepsilon })} (\tilde{\bm a},\tilde{\bm j}) \in 
\mathcal{A}_r \cdot E(\tilde{\bm a},\tilde{\bm j})$. \\ \\
(b) For 
${\bm \varepsilon}, 
\tilde{{\bm \varepsilon}} \in 
\mathcal{X}_r({\bm a},{\bm j}),$ we have 
$\mathcal{A}_r \cdot B^{(\tilde{\bm \varepsilon })} ({\bm a},{\bm j}) \subseteq 
\mathcal{A}_r \cdot B^{({\bm \varepsilon })} ({\bm a},{\bm j})
\mbox{ if and only if } {\bm \varepsilon } \leq \tilde{\bm \varepsilon }$ 
and  
$\mathcal{A}_r \cdot B^{(\tilde{\bm \varepsilon })} ({\bm a},{\bm j}) = 
\mathcal{A}_r \cdot B^{({\bm \varepsilon })} ({\bm a},{\bm j})
\mbox{ if and only if } {\bm \varepsilon } = \tilde{\bm \varepsilon }.$ \\ \\ 
(c) The $k$-algebra $\mathcal{A}_r \cdot E({\bm a},{\bm j}) $ has the elements 
$B^{({\bm \theta })} ({\bm a},{\bm j})$ with 
${\bm \theta } \in \mathcal{X}_r({\bm a},{\bm j})$ as a $k$-basis. \\ 

\noindent {\itshape Proof.} 
First we claim that the elements $B^{({\bm \theta })} ({\bm a},{\bm j})$ for 
${\bm \theta} \in \mathcal{X}_r({\bm a},{\bm j})$ are linearly independent over $k$. 
Note that $\mathcal{X}_1(a_0,j_0)$ is equal to $\{ 0,1\}$ if $j_0 \neq 0$ when $p$ is odd, or 
if $a_0 \neq 1$ when $p=2$, and to $\{ 0\}$ otherwise.  In the former case 
$B^{(0 )} (a_0,j_0)$ and $B^{(1 )} (a_0,j_0)$ are linearly independent over $k$ by 
Lemma 3.3. Hence the claim holds for 
$r=1$. 
Suppose that $r \geq 2$ and $$\sum_{{\bm \theta} \in \mathcal{X}_r({\bm a},{\bm j})} 
\alpha_{\bm \theta}B^{({\bm \theta })} ({\bm a},{\bm j})=0,$$ where 
$\alpha_{\bm \theta} \in k$. If we write  
$({\bm a}', {\bm j}') =\big( (a_i,j_i)\big)_{i=1}^{r-1} \in \mathcal{P}^{r-1}$, 
${\bm \theta} =(\theta_0, \dots, \theta_{r-1}) \in 
\mathcal{X}_r({\bm a},{\bm j})$ and  
${\bm \theta}' =(\theta_1, \dots, \theta_{r-1}) \in 
\mathcal{X}_{r-1}({\bm a}',{\bm j}')$, we have 
\begin{eqnarray*} 
0 &=& \sum_{{\bm \theta} \in \mathcal{X}_r({\bm a},{\bm j})} 
\alpha_{\bm \theta}B^{({\bm \theta })} ({\bm a},{\bm j}) 
= \sum_{{\bm \theta} \in \mathcal{X}_r({\bm a},{\bm j})} \alpha_{\bm \theta} 
Z^{(\theta_0)}\Big( B^{({\bm \theta}')}({\bm a}', {\bm j}') ; (a_0,j_0) \Big) \\
& = & \sum_{\theta_0 \in \mathcal{X}_1(a_0,j_0)} Z^{(\theta_0)}
\Big( \sum_{{\bm \theta}' \in \mathcal{X}_{r-1}({\bm a}',{\bm j}')} 
\alpha_{(\theta_0, {\bm \theta}')}
B^{({\bm \theta}')}({\bm a}', {\bm j}') ; (a_0,j_0) \Big),
\end{eqnarray*} 
where $(\theta_0, {\bm \theta}')$ means ${\bm \theta}$. 
By (a) in the remark of Lemma 3.4 we have 
$$Z^{(\theta_0)}
\Big( \sum_{{\bm \theta}' \in \mathcal{X}_{r-1}({\bm a}', {\bm j}')} 
\alpha_{(\theta_0, {\bm \theta}')}
B^{({\bm \theta}')}({\bm a}', {\bm j}') ; (a_0,j_0) \Big) =0,$$
and hence $\sum_{{\bm \theta}' \in \mathcal{X}_{r-1}({\bm a}', {\bm j}')} 
\alpha_{(\theta_0, {\bm \theta}')}
B^{({\bm \theta}')}({\bm a}', {\bm j}')=0$ 
for each $\theta_0 \in \mathcal{X}_1(a_0,j_0)$. 
Since $B^{({\bm \theta}')}({\bm a}', {\bm j}')$ with 
${\bm \theta}' \in \mathcal{X}_{r-1}({\bm a}',{\bm j}')$ are linearly independent by induction, 
we obtain $\alpha_{(\theta_0, {\bm \theta}')}=0$ for each 
${\bm \theta}' \in \mathcal{X}_{r-1}({\bm a}',{\bm j}')$. It follows that 
 $\alpha_{{\bm \theta}}=0$ for each 
${\bm \theta} \in \mathcal{X}_{r}({\bm a},{\bm j})$, 
and the claim follows.  

Next we claim that $\mathcal{A}_r \cdot B^{(\bm{\varepsilon})}(\bm{a}, \bm{j})$ 
is spanned by all $B^{(\bm{\theta})}(\bm{a}, \bm{j})$ with 
$\bm{\theta} \in \mathcal{X}_{r}({\bm a},{\bm j})$ and 
$\bm{\theta} \geq \bm{\varepsilon}$. Let $V$ be the subspace spanned by 
all $B^{(\bm{\theta})}(\bm{a}, \bm{j})$ with 
$\bm{\theta} \in \mathcal{X}_{r}({\bm a},{\bm j})$ and 
$\bm{\theta} \geq \bm{\varepsilon}$. Suppose that an element  
$\bm{\theta} \in \mathcal{X}_{r}({\bm a},{\bm j})$ satisfies 
$\bm{\theta} \geq \bm{\varepsilon}$. For an integer $s$ with $0 \leq s \leq r-1$, 
if $s$ satisfies $j_s \neq 0$ when $p$ is odd, or $a_s \neq 1$ (i.e. $(a_s,j_s)=(0,1/2)$) 
when $p=2$, and if $\theta_s =0$, then $\bm{\theta} + {\bm e}_{s+1} \in 
\mathcal{X}_{r}({\bm a},{\bm j})$ and $\bm{\theta} + {\bm e}_{s+1}  \geq \bm{\varepsilon}$. 
Thus we see that $Y^{(p^s)}X^{(p^s)} B^{(\bm{\theta})}(\bm{a}, \bm{j}) \in V$ 
by Proposition 3.7. Since $B^{(\bm{\theta})}(\bm{a}, \bm{j})$ is a 
$\mathcal{U}_r^0$-weight vector, $V$ is closed under the action of $\mathcal{A}_r$. 
Moreover, since $B^{(\bm{\varepsilon})}(\bm{a}, \bm{j}) \in V$, we obtain 
$\mathcal{A}_r \cdot B^{(\bm{\varepsilon})}(\bm{a}, \bm{j}) \subseteq V$. To show 
the reverse inclusion, we have to check that 
$B^{(\bm{\theta})}(\bm{a}, \bm{j}) \in \mathcal{A}_r \cdot B^{(\bm{\varepsilon})}(\bm{a}, \bm{j})$ for any $\bm{\theta} \in \mathcal{X}_{r}({\bm a},{\bm j})$ satisfying 
$\bm{\theta} \geq \bm{\varepsilon}$. It is clear when  
$d(\bm{\theta}, \bm{\varepsilon})=0$ (i.e. $\bm{\theta}=\bm{\varepsilon}$), so suppose that $d(\bm{\theta}, \bm{\varepsilon})>0$. There exists an integer $s$ with 
$0 \leq s \leq r-1$ 
such that $\varepsilon_s =0$ and $\theta_s=1$. For this integer $s$ note that 
$j_s \neq 0$ when $p$ is odd, or  
$a_s \neq 1$ (i.e. $(a_s, j_s)=(0,1/2)$) when $p=2$. Then 
$$Y^{(p^s)} X^{(p^s)} B^{({\bm \theta }-{\bm e}_{s+1})} ({\bm a},{\bm j}) = 
\Bigg( j_s^2 - \bigg( \dfrac{a_{s}+1}{2}\bigg)^2\Bigg) B^{({\bm \theta }-{\bm e}_{s+1})} 
({\bm a},{\bm j}) +4j_s^2 
B^{({\bm \theta })} ({\bm a},{\bm j})$$
by Proposition 3.7. Since ${\bm \theta }-{\bm e}_{s+1}$ is an element of 
$\mathcal{X}_{r}({\bm a},{\bm j})$ satisfying 
${\bm \theta }-{\bm e}_{s+1} \geq \bm{\varepsilon}$ and 
$d({\bm \theta }-{\bm e}_{s+1},\bm{\varepsilon})= d({\bm \theta }, \bm{\varepsilon})-1$, 
we obtain $B^{({\bm \theta }-{\bm e}_{s+1})} ({\bm a},{\bm j}) \in 
\mathcal{A}_r \cdot B^{(\bm{\varepsilon})}(\bm{a}, \bm{j})$ by induction. Moreover, since 
$4j_s^2 \neq 0$ in $\mathbb{F}_p$, we have 
$$B^{({\bm \theta })} ({\bm a},{\bm j}) = \dfrac{1}{4j_s^2} 
\Bigg( Y^{(p^s)} X^{(p^s)} -j_s^2 + \bigg( \dfrac{a_s+1}{2}\bigg)^2 \Bigg) 
B^{({\bm \theta }-{\bm e}_{s+1})} ({\bm a},{\bm j}) \in \mathcal{A}_r \cdot B^{(\bm{\varepsilon})}(\bm{a}, \bm{j}).$$
Therefore, we obtain $\mathcal{A}_r \cdot B^{(\bm{\varepsilon})}(\bm{a}, \bm{j})=V$ 
and the proof is complete. $\square$ \\

The following lemma enables us to determine radical series of the 
$\mathcal{A}_r$-modules  $\mathcal{A}_r \cdot  \nobreak B^{(\bm{\varepsilon})}(\bm{a}, \bm{j})$ 
with ${\bm \varepsilon} \in \mathcal{X}_r({\bm a},{\bm j})$. 

\begin{Lem}
Let 
${\bm \varepsilon} =(\varepsilon_0, \dots, \varepsilon_{r-1}), 
\tilde{{\bm \varepsilon}} =(\tilde{\varepsilon}_0, \dots, \tilde{\varepsilon}_{r-1}) 
\in \mathcal{X}_r({\bm a},{\bm j})$. Then the product 
$B^{({\bm \varepsilon })} ({\bm a},{\bm j})B^{(\tilde{\bm \varepsilon })} ({\bm a},{\bm j})$ 
is equal to zero if there is an integer $s$ with $0 \leq s \leq r-1$ such that  
$\varepsilon_s = \tilde{\varepsilon}_s =1$, and to 
$B^{({\bm \varepsilon }+\tilde{\bm \varepsilon })} ({\bm a},{\bm j})$ otherwise. 
In the latter case ${\bm \varepsilon }+\tilde{\bm \varepsilon }$ also lies in 
$\mathcal{X}_r({\bm a},{\bm j})$.
\end{Lem}

\noindent {\itshape Proof.} Suppose that there is an integer $s$ satisfying 
$\varepsilon_s = \tilde{\varepsilon}_s =1$. Note that $4j_s^2 \neq 0$ in 
$\mathbb{F}_p$. By Proposition 3.7 we have 
$$B^{({\bm \varepsilon })} ({\bm a},{\bm j}) = \dfrac{1}{4j_s^2} 
\Bigg( Y^{(p^s)} X^{(p^s)} -j_s^2 + \bigg( \dfrac{a_s+1}{2}\bigg)^2 \Bigg) 
B^{({\bm \varepsilon }-{\bm e}_{s+1})} ({\bm a},{\bm j})$$
and $Y^{(p^s)} X^{(p^s)}B^{(\tilde{\bm \varepsilon })} ({\bm a},{\bm j}) 
= \Big( j_s^2- \big((a_s+1)/2\big)^2\Big)B^{(\tilde{\bm \varepsilon })} ({\bm a},{\bm j}) $. 
Then 
\begin{eqnarray*}
\lefteqn{B^{({\bm \varepsilon })} ({\bm a},{\bm j}) B^{(\tilde{\bm \varepsilon })} ({\bm a},{\bm j})} 
\\ 
& = &  \dfrac{1}{4j_s^2} 
\Bigg( Y^{(p^s)} X^{(p^s)} -j_s^2 + \bigg( \dfrac{a_s+1}{2}\bigg)^2 \Bigg) 
B^{(\tilde{\bm \varepsilon })} ({\bm a},{\bm j})
B^{({\bm \varepsilon }-{\bm e}_{s+1})} ({\bm a},{\bm j}) \\
& = & 0.
\end{eqnarray*}

On the other hand, suppose that there are no integers $s$ satisfying 
$\varepsilon_s = \tilde{\varepsilon}_s =1$. Clearly 
${\bm \varepsilon }+\tilde{\bm \varepsilon }$ lies in 
$\mathcal{X}_r({\bm a},{\bm j})$ again. We prove the lemma by induction on 
$\mathcal{W}({\bm \varepsilon }+\tilde{\bm \varepsilon })$. It is clear when 
$\mathcal{W}({\bm \varepsilon }+\tilde{\bm \varepsilon })=0$, since 
${\bm \varepsilon } = \tilde{\bm \varepsilon }= (0,\dots,0)$. Suppose that 
$\mathcal{W}({\bm \varepsilon }+\tilde{\bm \varepsilon })>0$. We may assume that 
there is an integer $s$ such that $\varepsilon_s=1$. By induction, the product
$B^{({\bm \varepsilon }-{\bm e}_{s+1})} ({\bm a},{\bm j})
B^{(\tilde{\bm \varepsilon })} ({\bm a},{\bm j})$ is 
equal to $B^{({\bm \varepsilon }+\tilde{\bm \varepsilon }-{\bm e}_{s+1})} 
({\bm a},{\bm j})$. Thus we have 
\begin{eqnarray*}
\lefteqn{B^{({\bm \varepsilon })} ({\bm a},{\bm j}) B^{(\tilde{\bm \varepsilon })} ({\bm a},{\bm j})} 
\\ 
& = &  \dfrac{1}{4j_s^2} 
\Bigg( Y^{(p^s)} X^{(p^s)} -j_s^2 + \bigg( \dfrac{a_s+1}{2}\bigg)^2 \Bigg) 
B^{({\bm \varepsilon }-{\bm e}_{s+1})} ({\bm a},{\bm j}) 
B^{(\tilde{\bm \varepsilon })} ({\bm a},{\bm j})
\\
& = & \dfrac{1}{4j_s^2} 
\Bigg( Y^{(p^s)} X^{(p^s)} -j_s^2 + \bigg( \dfrac{a_s+1}{2}\bigg)^2 \Bigg) 
B^{({\bm \varepsilon }+\tilde{\bm \varepsilon }-{\bm e}_{s+1})} 
({\bm a},{\bm j})\\
& = & B^{({\bm \varepsilon }+\tilde{\bm \varepsilon })}({\bm a},{\bm j}), 
\end{eqnarray*}
as required. $\square$ \\

Since all $E({\bm a}, {\bm j})$ with $({\bm a}, {\bm j}) \in \mathcal{P}^r$ 
are central idempotents of $\mathcal{A}_r$ whose sum is 1, 
the representation theory for the algebras $\mathcal{A}_r \cdot E({\bm a}, {\bm j})$ 
completely determines that for  $\mathcal{A}_r$ 
(see \cite[ch. 1. Theorem 4.7]{nagao-tsushimabook}). For a fixed 
$({\bm a}, {\bm j}) \in \mathcal{P}^r$, 
set $w= \mathcal{W}({\bm j})$ if $p$ is 
odd, and $w= r-\mathcal{W}({\bm a})$ if $p=2$ (i.e. $w$ is the number of the integers $s$ 
with $0 \leq s \leq r-1$ satisfying $j_s \neq 0$ if $p$ is odd, or $a_s \neq 1$ if $p=2$). 
Then for ${\bm \varepsilon} \in \mathcal{X}_r({\bm a},{\bm j})$, 
the $\mathcal{A}_r$-module 
$\mathcal{A}_r \cdot B^{({\bm \varepsilon })} ({\bm a},{\bm j})$ has 
dimension $2^{w-\mathcal{W}({\bm \varepsilon})}$ by Theorem 3.8. In particular, 
the $k$-algebra $\mathcal{A}_r \cdot E({\bm a}, {\bm j})$ has dimension $2^w$ which is 
also the cardinality of $\mathcal{X}_r({\bm a},{\bm j})$. 

\begin{Prop}
For a positive integer $i$,  
we have 
\begin{eqnarray*}
\big( {\rm rad}(\mathcal{A}_r \cdot E({\bm a},{\bm j}))\big)^i
& = &\sum_{{\bm \theta } \in \mathcal{X}_r({\bm a},{\bm j}),\  
\mathcal{W}({\bm \theta })=i}
\mathcal{A}_r \cdot 
B^{({\bm \theta })} ({\bm a},{\bm j}) \\
& = & \sum_{{\bm \theta } \in \mathcal{X}_r({\bm a},{\bm j}),\  
\mathcal{W}({\bm \theta }) \geq i}
k \cdot 
B^{({\bm \theta })} ({\bm a},{\bm j}). 
\end{eqnarray*}
In particular, $\big( {\rm rad}(\mathcal{A}_r \cdot E({\bm a},{\bm j}))\big)^i=0$ if 
and only if 
$i>w$.
\end{Prop}

\noindent {\itshape Proof.} The second equality follows immediately from 
Theorem 3.8, so we only have to show the first equality. Lemma 3.9 implies that 
the subspace 
$$\sum_{{\bm \theta} \in \mathcal{X}_r({\bm a},{\bm j}),\  
\mathcal{W}({\bm \theta }) \geq 1}
k \cdot 
B^{({\bm \theta })} ({\bm a},{\bm j})$$
is an ideal of the algebra $\mathcal{A}_r \cdot E({\bm a},{\bm j})$ and that 
a product of 
 $w+1$ elements in the subspace 
is equal to 0. Hence the subspace  is a nilpotent 
ideal of $\mathcal{A}_r \cdot E ({\bm a},{\bm j})$. Moreover, by Theorem 3.8  
we see that the nilpotent ideal has codimension one in 
$\mathcal{A}_r \cdot E({\bm a},{\bm j})$ and hence is equal to 
${\rm rad}(\mathcal{A}_r \cdot E({\bm a},{\bm j}))$. Thus the result for $i=1$ follows. 
The result for arbitrary $i$ follows 
easily from Lemma 3.9 using induction on $i$.  $\square$ \\ 

Let 
${\bm \tau} = {\bm \tau}({\bm a},{\bm j})= (\tau_0, \dots, \tau_{r-1}) \in 
\mathcal{X}_r({\bm a},{\bm j})$ be the element 
such that $\tau_s =0$ if $j_s =0$ when $p$ is odd or  $a_s =1$ when $p=2$, 
and $\tau_s =1$ otherwise for $0 \leq s \leq r-1$. This is a unique element of 
$\mathcal{X}_r({\bm a},{\bm j})$ which has the largest Hamming weight $w$. 
Set $S_{({\bm a}, {\bm j})} = \mathcal{A}_r \cdot 
B^{({\bm \tau})}({\bm a}, {\bm j}) = k \cdot 
B^{({\bm \tau})}({\bm a}, {\bm j}) $. Clearly this is a simple $\mathcal{A}_r$-module.

Now we give radical series of the $\mathcal{A}_r$-modules  
$\mathcal{A}_r \cdot B^{({\bm \varepsilon })} ({\bm a},{\bm j})$ for 
${\bm \varepsilon } \in \mathcal{X}_r({\bm a}, {\bm j})$. 

\begin{The}
Let 
${\bm \varepsilon} \in \mathcal{X}_r({\bm a},{\bm j})$. For $i \in \mathbb{Z}_{\geq 0}$ 
we have 
\begin{eqnarray*}
{\rm rad}^i_{\mathcal{A}_r}(\mathcal{A}_r \cdot B^{({\bm \varepsilon })} ({\bm a},{\bm j}))
&=& \sum_{{\bm \theta } \in \mathcal{X}_r({\bm a}, {\bm j}),\ 
d({\bm \theta},{\bm \varepsilon})=i,\ {\bm \theta }\geq {\bm \varepsilon }}
\mathcal{A}_r \cdot 
B^{({\bm \theta })} ({\bm a},{\bm j}) \\
&=& \sum_{{\bm \theta } \in \mathcal{X}_r({\bm a}, {\bm j}),\ 
d({\bm \theta},{\bm \varepsilon}) \geq i,\ {\bm \theta }\geq {\bm \varepsilon }}
k \cdot 
B^{({\bm \theta })} ({\bm a},{\bm j}),
\end{eqnarray*}
and the Loewy length of 
$\mathcal{A}_r \cdot B^{({\bm \varepsilon })} ({\bm a},{\bm j})$ is 
$w+1-\mathcal{W}({\bm \varepsilon})$. Moreover, for an integer $i$ with 
$0 \leq i \leq w-\mathcal{W}({\bm \varepsilon})$ the quotient 
$$
{\rm rad}^{i}_{\mathcal{A}_r}(\mathcal{A}_r \cdot B^{({\bm \varepsilon })} ({\bm a},{\bm j}))/ 
{\rm rad}^{i+1}_{\mathcal{A}_r}(\mathcal{A}_r \cdot B^{({\bm \varepsilon })} ({\bm a},{\bm j})) 
$$
is isomorphic to a direct sum of ${w-\mathcal{W}({\bm \varepsilon}) \choose i} $ copies of 
$S_{({\bm a},{\bm j})}$.
In particular, $\mathcal{A}_r \cdot B^{({\bm \varepsilon })} ({\bm a},{\bm j})$ is an 
indecomposable $\mathcal{A}_r$-module whose head is isomorphic to 
$S_{({\bm a},{\bm j})}$.
\end{The}

\noindent {\itshape Proof.} The first statement follows  easily from Proposition 3.10 
and Lemma 3.9 since 
$${\rm rad}^{i}_{\mathcal{A}_r}(\mathcal{A}_r \cdot 
B^{({\bm \varepsilon })} ({\bm a},{\bm j}))= 
\big( {\rm rad}(\mathcal{A}_r \cdot E({\bm a},{\bm j}))\big)^i \cdot 
B^{({\bm \varepsilon })} ({\bm a},{\bm j})$$
if $i \geq 1$. Then the  second statement follows easily from the first statement and Proposition 
3.7. $\square$ \\

This theorem implies that 
all the idempotents 
$E({\bm a}, {\bm j})$ are primitive, and all $S_{({\bm a}, {\bm j})}$ form a complete set of 
non-isomorphic simple $\mathcal{A}_r$-modules. In particular, 
$\mathcal{A}_r \cdot E({\bm a}, {\bm j})$ is a block algebra of $\mathcal{A}_r$ which has 
$S_{({\bm a}, {\bm j})}$ as a unique simple 
$\mathcal{A}_r \cdot E({\bm a}, {\bm j})$-module. 

The following theorem shows that the $\mathcal{A}_r$-modules  
$\mathcal{A}_r \cdot B^{({\bm \varepsilon })} ({\bm a},{\bm j})$ 
for ${\bm \varepsilon} \in \mathcal{X}_r({\bm a},{\bm j})$ are rigid 
(i.e. have identical radical and socle series). 

\begin{The}
Let 
${\bm \varepsilon} \in \mathcal{X}_r({\bm a},{\bm j})$. For an integer $i$ with 
$0 \leq i \leq  w - \mathcal{W}({\bm \varepsilon})$ 
we have 
$${\rm rad}^i_{\mathcal{A}_r}(\mathcal{A}_r \cdot B^{({\bm \varepsilon })} 
({\bm a},{\bm j})) = 
{\rm soc}^{w+1-\mathcal{W}({\bm \varepsilon})-i}_{\mathcal{A}_r}
(\mathcal{A}_r \cdot B^{({\bm \varepsilon })} ({\bm a},{\bm j})).$$
In particular, the socle of the 
$\mathcal{A}_r$-module  
$\mathcal{A}_r \cdot B^{({\bm \varepsilon })} ({\bm a},{\bm j})$ is 
isomorphic to $S_{({\bm a}, {\bm j})}$. 
\end{The}

\noindent {\itshape Proof.} It is clear when $i=0$ since the $\mathcal{A}_r$-module 
$\mathcal{A}_r \cdot B^{({\bm \varepsilon })} ({\bm a},{\bm j})$ has Loewy length 
$w+1-\mathcal{W}({\bm \varepsilon})$. So we may assume $i \geq 1$. It is enough to 
show that $${\rm soc}^{w+1-\mathcal{W}({\bm \varepsilon})-i}_{\mathcal{A}_r}
(\mathcal{A}_r \cdot B^{({\bm \varepsilon })} ({\bm a},{\bm j})) 
\subseteq {\rm rad}^i_{\mathcal{A}_r}(\mathcal{A}_r \cdot B^{({\bm \varepsilon })} 
({\bm a},{\bm j})). $$
Suppose that an element $u \in \mathcal{A}_r \cdot B^{({\bm \varepsilon })} 
({\bm a},{\bm j})$ does not lie in ${\rm rad}^i_{\mathcal{A}_r}(\mathcal{A}_r \cdot B^{({\bm \varepsilon })} 
({\bm a},{\bm j}))$. We only have to check that 
$u \not\in {\rm soc}^{ w+1-\mathcal{W}({\bm \varepsilon })-i}_{\mathcal{A}_r}
(\mathcal{A}_r \cdot 
B^{({\bm \varepsilon})} ({\bm a},{\bm j}))$. By Theorem 3.8, $u$ can be written as a $k$-linear combination 
of the elements $B^{({\bm \theta })} ({\bm a},{\bm j})$ with 
${\bm \theta } \in \mathcal{X}_r({\bm a},{\bm j})$ and 
${\bm \theta } \geq {\bm \varepsilon }$. By the assumption of $u$ and by Theorem 
3.11,  if we  choose $\tilde{\bm \theta } \in \mathcal{X}_r({\bm a},{\bm j})$  where 
the coefficient of $B^{(\tilde{\bm \theta })} ({\bm a},{\bm j})$ in $u$ is nonzero 
such that $d(\tilde{\bm \theta }, {\bm \varepsilon })$ is minimal, 
$d(\tilde{\bm \theta }, {\bm \varepsilon })$ must be smaller than $i$, and hence 
$\mathcal{W}(\tilde{\bm \theta })-\mathcal{W}({\bm \varepsilon })  \leq i-1$.  Then 
$B^{({\bm \tau}-\tilde{\bm \theta })} ({\bm a},{\bm j}) u$ is a nonzero multiple of 
$B^{({\bm \tau})} ({\bm a},{\bm j})$ by Lemma 3.9. Since 
${\bm \tau}-\tilde{\bm \theta } \in \mathcal{X}_r({\bm a},{\bm j})$ and 
$$\mathcal{W}({\bm \tau}-\tilde{\bm \theta }) = 
\mathcal{W}({\bm \tau })-\mathcal{W}(\tilde{\bm \theta }) =
w-\mathcal{W}(\tilde{\bm \theta }) \geq w+1-\mathcal{W}({\bm \varepsilon })-i,$$
the element $B^{({\bm \tau}-\tilde{\bm \theta })} ({\bm a},{\bm j})$ must lie in 
$$\sum_{{\bm \theta } \in \mathcal{X}_r({\bm a}, {\bm j}),\ 
\mathcal{W}({\bm \theta }) \geq w+1-\mathcal{W}({\bm \varepsilon })-i}
k \cdot 
B^{({\bm \theta })} ({\bm a},{\bm j})
=\big( {\rm rad}(\mathcal{A}_r \cdot 
E ({\bm a},{\bm j})) \big)^{ w+1-\mathcal{W}({\bm \varepsilon })-i}.$$
This means $u \not\in {\rm soc}^{ w+1-\mathcal{W}({\bm \varepsilon })-i}_{\mathcal{A}_r}(\mathcal{A}_r \cdot 
B^{({\bm \varepsilon})} ({\bm a},{\bm j}))$. Therefore, the result follows. $\square$ \\

Actually, the $k$-algebra $\mathcal{A}_r \cdot E({\bm a}, {\bm j})$ is symmetric.

\begin{The}
$\mathcal{A}_r \cdot E({\bm a}, {\bm j})$ is a symmetric $k$-algebra. 
\end{The}

\noindent{\itshape Proof.} 
 Let $f : \mathcal{A}_r \cdot E({\bm a}, {\bm j}) 
\rightarrow k$ be the  $k$-linear map defined as follows: 
for  ${\bm \varepsilon} 
\in \mathcal{X}_r({\bm a},{\bm j})$, $f(B^{({\bm \varepsilon })} ({\bm a},{\bm j})) =0$ 
if ${\bm  \varepsilon} \neq {\bm \tau}$, 
 and $f(B^{({\bm \tau })} ({\bm a},{\bm j})) =1$.   Let $u$ be a nonzero element of $\mathcal{A}_r \cdot E ({\bm a},{\bm j})$. 
By Theorem 3.8, $u$ can be written as 
$\sum_{{\bm \varepsilon} \in \mathcal{X}_r({\bm a},{\bm j})}\alpha_{\bm \varepsilon}B^{({\bm \varepsilon })} ({\bm a},{\bm j})$, $\alpha_{\bm \varepsilon} \in k$.  
Choose an element ${\bm \theta} \in \mathcal{X}_r({\bm a},{\bm j})$ 
with $\alpha_{\bm \theta} \neq 0$ such that 
$\mathcal{W}({\bm \theta})$ is minimal. 
Since ${\bm \tau }- {\bm \theta} \in \mathcal{X}_r({\bm a},{\bm j})$, the element 
$B^{({\bm \tau }- {\bm \theta})} ({\bm a},{\bm j})$ lies in 
$\mathcal{A}_r \cdot E({\bm a}, {\bm j}) $. 
By Lemma 3.9, 
we see that $B^{({\bm \tau }- {\bm \theta})} ({\bm a},{\bm j})u = 
\alpha_{\bm \theta} B^{({\bm \tau })} ({\bm a},{\bm j})$. 
This fact implies that ${\rm Ker} f$ 
 contains no nonzero ideals of $\mathcal{A}_r \cdot E({\bm a}, {\bm j})$. 
Thus $\mathcal{A}_r \cdot E({\bm a}, {\bm j})$ is a Frobenius algebra 
(see \cite[ch. 2. Theorem 8.13]{nagao-tsushimabook}). But it is 
also symmetric since 
$\mathcal{A}_r$ (hence $\mathcal{A}_r \cdot E({\bm a}, {\bm j})$) is commutative. 
 $\square$

\end{document}